\documentclass[10pt]{article} 

\usepackage{latexsym,amsmath,amssymb,a4wide,amscd}
\usepackage[latin1]{inputenc}
\usepackage[T1]{fontenc}

\usepackage{array}

\usepackage{enumerate}

\usepackage{amsthm}

\usepackage{lscape}

\usepackage{url}

\usepackage{mathabx}

\usepackage{longtable}
\usepackage{multirow}

\usepackage{MnSymbol}

\usepackage{amsthm}
\numberwithin{equation}{section}
\numberwithin{figure}{section}
\theoremstyle{plain}
\newtheorem {satz}{Theorem}[section]

\newtheorem {prop}[satz]{Proposition}

\newtheorem {conj}[satz]{Conjecture}

\theoremstyle{definition}

\theoremstyle{remark}

\usepackage[pdftex]{graphicx} 

\usepackage{pdfpages}

\usepackage[nooneline]{caption}

\title{Computationally proving triangulated 4-manifolds to be diffeomorphic}

\author{Benjamin A.~Burton
        \and
        Jonathan Spreer\thanks{%
            School of Mathematics and Physics,
            The University of Queensland,
            AUS.
            \texttt{bab@maths.uq.edu.au}, \texttt{j.spreer@uq.edu.au}}}

\begin{document}
\maketitle

\begin{abstract}
We present new computational methods for proving diffeomorphy of
triangulated 4-manifolds,
including algorithms and topological software that
can for the first time effectively handle the complexities that arise in
dimension four and be used for large scale experiments.
\end{abstract}

{\small \textbf{MSC 2000: } Primary 57Q25; 
Secondary 57Q15, 
14J28, 
68T20 

}

\section{Introduction}

	In dimensions $\leq 3$, every topological manifold
	has a unique smooth structure up to diffeomorphism.
	In dimensions $\geq 4$ this is no longer true:
	there are pairs of $4$-manifolds which are {\em homeomorphic}
	(they represent the same topological manifold) but not {\em diffeomorphic}
	(they represent two distinct smooth manifolds) 
	\cite{Gompf}. Finding such pairs is important; indeed,
	the only outstanding variant of the Poincar{\'e} conjecture
	asks whether one can find two
	non-diffeomorphic smooth 4-spheres \cite{Gompf}.

	We move this problem to the \emph{piecewise linear}
	setting, which is better suited for computation.
	Here manifolds are given as {\em triangulations} (decompositions
	into simplices).
	Piecewise linear manifolds are in 1-to-1 correspondence with smooth
	manifolds for dimensions $\leq 6$, and so results translate
	between both settings; we use both languages interchangably in this paper.

	Despite this equivalence, work on non-diffeomorphic pairs is done
	exclusively in the smooth setting. The only example 
	of two 4-manifold \emph{triangulations} that are
	homeomorphic but not diffeomorphic follows a
	well-established result for the smooth setting \cite{Benedetti13RDMT}.
	One of the few candidates from the PL-world is a pair of triangulations
	of the {\em K3 surface} (one of the four fundamental building
	blocks of simply connected 4-manifolds):
	the 16-vertex $(K3)_{16}$ of Casella and K{\"u}hnel
	\cite{Casella01TrigK3MinNumVert},
	and the 17-vertex
	$(K3)_{17}$ of Spreer and K{\"u}hnel \cite{Spreer09CombPorpsOfK3}.
	The smooth type of $(K3)_{17}$ is canonical, but the smooth type
	of $(K3)_{16}$ remains unknown.  It is conjectured
	\cite{Spreer09CombPorpsOfK3}:

	\vspace{-2mm}
	\begin{conj}
		\label{conj:main}
		$(K3)_{16}$ and $(K3)_{17}$ are diffeomorphic.
	\end{conj}
	
	In the computational setting,
	proving that $(K3)_{16}$ and $(K3)_{17}$ are
	\emph{homeomorphic} is easy, using the software
	\texttt{simpcomp} \cite{simpcompISSAC} in conjunction with
	Freedman's celebrated classification of simply connected
	$4$-manifolds \cite{Freedman82Top4DimMnf}.
	Proving they are \emph{diffeomorphic}
	is much harder. Our approach is based on a theorem of Pachner
	\cite{Pachner87KonstrMethKombHomeo}, which states
	that two triangulated manifolds are diffeomorphic if and only if they
	are related by a sequence of {\em bistellar moves} (local
	modifications).

	Bistellar moves offer significant challenges in dimensions $\geq 4$:
	``effective'' sequences of moves can be extremely difficult to find.
	Indeed, the \emph{number} of moves required to connect two diffeomorphic
	triangulations of size $n$ must have no computable upper bound 
	\cite{Mijatovic03SimplTrigsOfS3}.
	
	Here we describe work in progress towards resolving
	Conjecture \ref{conj:main}, including effective heuristics and fast
	algorithms for manipulating 4-manifold triangulations with bistellar flips,
	and a tight lower bound on the size of a 1-vertex triangulation
	of the $K3$ surface with over a million distinct realisations.
	This work has wider relevance, and forms the beginning of a larger
	project to explicitly construct and study ``exotic'' triangulated
	4-manifolds.


	\section{Minimal triangulations of the K3 surface}

	For computation, we want to
	triangulate manifolds using few top-dimensional simplices.
	We therefore work with
	\emph{generalised triangulations}, which are collections
	of abstract simplices whose facets are identified in pairs---these
	can allow far fewer simplices than the more rigid simplicial complexes.
	We also favour triangulations with just one vertex, which in lower
	dimensions offer significant advantages for both theory and computation.

	\emph{Proving} minimality is extremely difficult in three dimensions.
	In four dimensions, we solve this completely for the $K3$ surface
	in the one-vertex setting:

	\vspace{-2mm}
	\begin{prop}
		\label{prop:lowerBound}
		For any triangulation of the $K3$ surface we have
		$f_4 \geq 146 -6f_0$,
		where $f_4$ denotes the number of $4$-dimensional simplices
		and $f_0$ denotes the number of vertices.
		In the case where $f_0=1$, this bound is tight.
	\end{prop}
	
	We prove $f_4 \geq 146 -6f_0$ by combining (i)~the fact that the $K3$
	surface is simply connected and has Euler characteristic
	$\chi (K3) = f_0 - f_1 + f_2 - f_3 + f_4 = 24$,
	with (ii)~the {\em Dehn-Sommer\-ville equations}
	$2f_1 - 3f_2 + 4f_3 - 5f_4 = 0$ and
	$2 f_3 - 5 f_4 = 0$.
	Here each $f_i$ denotes the number of $i$-faces of the triangulation.


	We prove this bound is tight by reducing both
	$(K3)_{16}$ and $(K3)_{17}$ using bistellar moves to
	one-vertex triangulations with $(f_0,f_1,f_2,f_3,f_4)=(1,1,234,350,140)$.
	In three dimensions, such a ``simplification'' of triangulations
	is fast and effective \cite{Burton12CompTopWRegina},
	but in four dimensions it is far more difficult
	and requires the interaction of many different tools and heuristics.

	Our approach incorporates:
	(i)~classical greedy techniques, which reduce a triangulation
	as far as possible using local moves;
	(ii)~``composite'' moves that collapse edges and
	reduce triangulations near low-degree edges and triangles;
	(iii)~simulated annealing \cite{Bjoerner00SimplMnfBistellarFlips},
	where we apply the inverses of reducing moves
	to escape local minima;
	(iv)~breadth-first searching through the Pachner graph (or ``flip graph'')
	\cite{Burton12CompTopWRegina}.

	We note that the \emph{interaction} between these
	techniques is crucial: each technique failed to reduce
	$(K3)_{16}$ and $(K3)_{17}$ on its own.
	Again we contrast this with three dimensions, where these
	techniques are found to be highly effective even in isolation.

	All computations were performed using the new $4$-manifold toolkit
	in the software package \texttt{Regina} \cite{Burton09Regina}.

	\section{Connecting triangulations $(K3)_{16}$ and $(K3)_{17}$}

	To prove Conjecture \ref{conj:main} we must find a sequence of 
	local modifications connecting $(K3)_{16}$ and $(K3)_{17}$.
	In three dimensions, the following approach is often successful:
	(i)~simplify both triangulations as far as possible, and then
	(ii)~repeatedly apply random local modifications that preserve
	the number of simplices until both triangulations are identical.
	This often succeeds (provided both triangulations represent
	the same manifold) because many manifolds appear to have only
	few distinct minimal triangulations.

	In contrast, for the $K3$ surface the number of minimal triangulations 
	appears to be \emph{much} larger, and so the classical approach above
	does not work.  Instead we use a more sophisticated method
	to ensure that (i) every minimal triangulation 
	is visited only once and (ii) we can detect if a longer ``detour''
	through larger triangulations is required.
	Specifically, we run a dual-source breadth-first search through the
	{\em Pachner graph}, whose nodes represent triangulations of the $K3$ surface
	and whose arcs represent
	bistellar flips that preserve the number of simplices.
	The two sources are our minimal one-vertex triangulations
	of $(K3)_{16}$ and $(K3)_{17}$.

	Each time we perform a local move, we must test whether the resulting
	triangulation has been seen before (up to combinatorial isomorphism).
	For this we compute the
	{\em isomorphism signature} of the triangulation \cite{Burton11FlipGraph},
	a polynomial-time computable hash that uniquely identifies
	the isomorphism type.
	This reduces the comparison to a fast lookup,
	and the overall algorithm runs in time
	$O(T \log T \cdot n^2 \log n)$,
	where $T$ is the (large) number of triangulations,
	and $n$ is the (small) number of simplices in each.
	The search parallelises well, since the bottlenecks are the
	hashing and performing local moves.

	Thus far, the algorithm has detected
	$1 \, 738 \, 260$ distinct minimal one-vertex triangulations of
	the $K3$ surface.  The search is ongoing,
	and has neither exhausted the list of minimal
	triangulations nor connected $(K3)_{16}$ with $(K3)_{17}$.
	This enormous number of minimal triangulations is both interesting
	and surprising, offering a stark contrast to observations
	from dimension three.

	\section{Conclusion and future research}

	Proving Conjecture \ref{conj:main} would
	eliminate an important candidate for a pair of homeomorphic
	but non-diffeomorphic simply connected $4$-manifolds.
	Moreover, as noted earlier, this work has a wider appeal:
	it shows for the first time how difficult problems of
	diffeomorphism and ``exotic structures'' in 4-manifold
	topology can be realistically tackled using computational tools.

	Future developments will include: multiple-vertex
	triangulations containing fewer $4$-simplices;
	a richer set of local modifications; and
	distributed algorithms for use on high-performance computing facilities.



\small
\bibliographystyle{abbrv}

\end{document}